\documentclass[12pt]{amsart}

\usepackage{fullpage}
\usepackage{amssymb}

\newtheorem{theorem}{Theorem}
\newtheorem{lemma}{Lemma}[section]

\theoremstyle{remark}

\newtheorem*{notation}{Notation}

\numberwithin{equation}{section}

\newcommand{\eps}{\varepsilon}

\title[Weyl sums over primes]
      {On Weyl sums over primes in short intervals}
\author{Angel V. Kumchev}
\address{Department of Mathematics,
         Towson University,
         7800 York Road,
         Towson, MD 21252}
\email{akumchev@towson.edu}

\subjclass[2010]{Primary 11L07; Secondary 11L15 11L20 11N36.}

\date{Draft from \today.}

\begin{document}

\maketitle

\section{Introduction}

In this note we pursue bounds for exponential sums of the form
\begin{equation}\label{i.1}
  f_k(\alpha; x, y) = \sum_{x < n \le x + y} \Lambda(n)e\left( \alpha n^k \right), 
\end{equation}
where $k \ge 2$ is an integer, $2 \le y \le x$, $\Lambda(n)$ is von Mangoldt's function, and $e(z) = e^{2\pi iz}$. When $y = x^{\theta}$ with $\theta < 1$, such exponential sums play a central role in applications of the Hardy--Littlewood circle method to additive problems with almost equal prime unknowns (see \cite{LvTa10, LvXu07, SuTa09}). When $\alpha$ is closely approximated by a rational number with a small denominator (i.e., when $\alpha$ is on a ``major arc''), Liu, L\"u and Zhan \cite{LiLvZh06} bounded $f_k(\alpha; x, x^{\theta})$ using methods from multiplicative number theory. Their result, which generalizes earlier work by Ren \cite{Ren05}, can be stated as follows.

\begin{theorem}\label{th1}
  Let $k \ge 1$, $7/10 < \theta \le 1$ and $0 < \rho \le \min\{ (8\theta - 5)/(6k+6), (10\theta-7)/15 \}$. Suppose that $\alpha$ is real and that there exist integers $a$ and $q$ satisfying
  \begin{equation}\label{i.2}
    1 \le q \le P, \quad (a,q) = 1, \quad |q\alpha - a| \le x^{-k + 2(1-\theta)}P,
  \end{equation} 
  with $P = x^{2k\rho}$. Then, for any fixed $\eps > 0$, 
  \[
    f_k\left( \alpha; x, x^{\theta} \right) \ll x^{\theta-\rho+\eps} + x^{\theta+\eps}\Xi(\alpha)^{-1/2},
  \]
  where $\Xi(\alpha) = q + x^{k-2(1-\theta)}|q\alpha - a|$.
\end{theorem}

For a given $P$, let $\mathfrak M(P)$ denote the set of real $\alpha$ that have rational approximations of the form \eqref{i.2}, and let $\mathfrak m(P)$ denote the complement of $\mathfrak M(P)$. In the terminology of the circle method, $\mathfrak M(P)$ is a set of \emph{major arcs} and $\mathfrak m(P)$ is the respective set of \emph{minor arcs}. The main goal of this note is to bound $f_k(\alpha; x, x^{\theta})$, $k \ge 3$, on sets of minor arcs by extending a theorem of the author \cite[Theorem 1]{Kumc06}, which gives the best known bound for $f_k(\alpha; x, x)$. We first state our result for cubic sums.

\begin{theorem}\label{th2}
  Let $\theta$ be a real number with $4/5 < \theta \le 1$ and suppose that $0 < \rho \le \rho_3(\theta)$, where
  \[
    \rho_3(\theta) = \min\left(\textstyle \frac 1{14}(2\theta-1), \frac 1{30}(14\theta-11), \frac 16(5\theta-4) \right).
  \]
  Then, for any fixed $\eps > 0$,
  \begin{equation}\label{i.3}
    \sup_{\alpha \in \mathfrak m(P)} \left|f_3 \left( \alpha; x, x^{\theta} \right)\right| 
    \ll x^{\theta - \rho + \eps} + x^{\theta+\eps}P^{-1/2}. 
  \end{equation}
\end{theorem}

We remark that when $\theta = 1$, Theorem \ref{th2} recovers the bound
\[
  \sup_{\alpha \in \mathfrak m(P)} \left|f_3(\alpha; x, x)\right| 
  \ll x^{13/14 + \eps} + x^{1+\eps}P^{-1/2},
\]
which is the essence of the cubic case of \cite[Theorem 3]{Kumc06}. In the case $k \ge 4$, our estimates take the following form.

\begin{theorem}\label{th3}
  Let $k \ge 4$ be an integer and $\theta$ be a real number with $1 - (k+2)^{-1} < \theta \le 1$. Suppose that $0 < \rho \le \rho_k(\theta)$, where
  \[
    \rho_k(\theta) = \min\left(\textstyle \frac 16\sigma_k(3\theta - 1), \frac 16((k+2)\theta-(k+1)) \right),
  \]
  with $\sigma_k$ defined by $\sigma_k^{-1} = \min( 2^{k-1}, 2k(k-2) )$. Then, for any fixed $\eps > 0$,
  \begin{equation}\label{i.4}
    \sup_{\alpha \in \mathfrak m(P)} \left|f_k \left( \alpha; x, x^{\theta} \right)\right| 
    \ll x^{\theta - \rho + \eps} + x^{\theta+\eps}P^{-1/2}. 
  \end{equation}
\end{theorem}

When $\theta = 1$ and $k \le 7$, this theorem also recovers the respective cases of \cite[Theorem~3]{Kumc06}. On the other hand, when $k \ge 8$, \eqref{i.4} is technically new even in the case $\theta = 1$, as we use the occasion to put on the record an almost automatic improvement of the theorems in \cite{Kumc06} that results from a recent breakthrough by Wooley \cite{Wool12a, Wool12b}.  

\begin{notation}
  Throughout the paper, the letter $\varepsilon$ denotes a sufficiently small positive real number. Any statement in which $\varepsilon$ occurs holds for each positive $\varepsilon$, and any implied constant in such a statement is allowed to depend on $\varepsilon$. The letter $p$, with or without subscripts, is reserved for prime numbers. As usual in number theory, $\mu(n)$, $\tau(n)$ and $\|x\|$ denote, respectively, the M\"obius function, the number of divisors function and the distance from $x$ to the nearest integer. We write $(a,b)=\gcd(a,b)$, and we use $m \sim M$ as an abbreviation for the condition $M < m \le 2M$.  
\end{notation}

\section{Auxiliary results}
\label{s2}

When $k \ge 3$, we define the multiplicative function $w_k(q)$ by
\[
  w_k \left( p^{ku + v} \right) =
  \begin{cases}
    kp^{-u - 1/2}, & \text{if } u \ge 0, v = 1, \\
    p^{-u - 1},    & \text{if } u \ge 0, v = 2, \dots, k.
  \end{cases}
\]
By the argument of \cite[Theorem 4.2]{Vaug97}, we have
\begin{equation}\label{ii.1}
  \sum_{1 \le x \le q} e\left( ax^k/q \right) \ll qw_k(q) \ll q^{1-1/k}
\end{equation}
whenever $k \ge 2$ and $(a, q) = 1$. We also need several estimates for sums involving the function $w_k(q)$. We list those in the following lemma.

\begin{lemma}\label{l1}
  Let $w_k(q)$ be the multiplicative function defined above. Then the following inequalities hold for any fixed $\eps > 0$:
  \begin{gather}
    \sum_{q \sim Q} w_k(q)^j \ll
    \begin{cases}
      Q^{-1 + \eps} & \text{if } k = 3, j = 4, \\
      Q^{-1 + 1/k}  & \text{if } k \ge 4, j = k;
    \end{cases}\label{ii.2}\\
    \sum_{n \sim N} w_k \left( \frac q {(q, n^j)} \right)
    \ll q^{\eps}w_k(q) N \qquad (1 \le j \le k); \label{ii.3}\\
    \sum_{ \substack{ n \sim N\\ (n, h) = 1}} w_k \left( \frac q{(q, R(n, h))} \right) 
    \ll q^{\eps}w_k(q)N + q^{\eps}, \label{ii.4}
  \end{gather}
  where $R(n,h) = \left( (n+h)^k - n^k \right)/h$.  
\end{lemma}

\begin{proof}
  See Lemmas 2.3 and 2.4 and inequality (3.11) in Kawada and Wooley \cite{KaWo01}.
\end{proof}

\begin{lemma}\label{l2}
  Let $k \ge 3$ be an integer and let $0 < \rho \le \sigma_k$, where $\sigma_k^{-1} = \min( 2^{k-1}, 2k(k-2) )$. Suppose that $y \le x$ and $x^k \le y^{k+1-2\rho}$. Then either
  \begin{equation}\label{ii.5}
    \sum_{x < n \le x+y} e \left( \alpha n^k \right) \ll y^{1-\rho + \eps},
  \end{equation}
  or there exist integers $a$ and $q$ such that
  \begin{equation}\label{ii.6} 
    1 \le q \le y^{k\rho}, \quad (a, q) = 1, \quad |q\alpha - a| \le x^{1-k}y^{k\rho - 1},
  \end{equation}
  and
  \begin{equation}\label{ii.7}
    \sum_{x < n \le x + y} e \left( \alpha n^k \right)
    \ll  \frac {w_k(q)y}{1 + yx^{k-1}|\alpha - a/q|} + x^{k/2 + \eps}y^{(1-k)/2}.
  \end{equation}
\end{lemma}

\begin{proof}
  By Dirichlet's theorem on Diophantine approximation, there exist integers $a$ and $q$ with
  \begin{equation}\label{ii.8}
    1 \le q \le y^{k - 1}, \quad (a, q) = 1, \quad |q\alpha - a| \le y^{1 - k}.
  \end{equation}
  When $q > y$, we rewrite the sum on the left of \eqref{ii.5} as
  \[
    \sum_{1 \le n \le y} e\left( \alpha n^k + \alpha_{k-1}n^{k-1} + \dots + \alpha_0 \right),
  \]
  where $\alpha_j = \binom kj\alpha [x]^{j-1}$. Hence, \eqref{ii.5} follows from Weyl's bound
  \[
    \sum_{1 \le n \le y} e\left( \alpha n^k + \alpha_{k-1}n^{k-1} + \dots + \alpha_0 \right) \ll y^{1-\sigma_k + \eps}.
  \]
  Under \eqref{ii.8}, this follows from \cite[Lemma 2.4]{Vaug97} when $\sigma_k = 2^{1-k}$ and from Wooley's recent improvement \cite{Wool12b} of Vinogradov's mean-value theorem otherwise. When $q \le X$, we deduce \eqref{ii.7} from \cite[Lemmas 6.1 and 6.2]{Vaug97} and \eqref{ii.1}. Thus, at least one of \eqref{ii.5} and \eqref{ii.7} holds. The lemma follows on noting that when conditions \eqref{ii.6} fail, inequality \eqref{ii.5} follows from \eqref{ii.7} and the hypothesis $x^k \le y^{k+1-2\rho}$.
\end{proof}

The following lemma is a slight variation of \cite[Lemma 6]{BaHa91}. The proof is the same.

\begin{lemma}\label{l3}
  Let $q$ and $N$ be positive integers exceeding $1$ and let $0 < \delta < \frac 12$. Suppose that $q \nmid a$ and denote by $\mathcal S$ the number of integers $n$ such that
  \[
    N < n \le 2N, \quad (n, q) = 1, \quad \left\| an^k/q \right\| < \delta.
  \]
  Then
  \[
    \mathcal S \ll \delta q^{\eps} (q + N).
  \]
\end{lemma}

\section{Multilinear Weyl sums}
\label{s3}

We write
\[
  \delta = x^{\theta-1}, \quad L = \log x, \quad \mathcal I = \left( x, x+x^{\theta} \right].
\]
We also set
\begin{equation}\label{iii.1}
  Q = \left( \delta x^{k-2\rho} \right)^{k/(2k-1)}.
\end{equation}
Recall that, by Dirichlet's theorem on Diophantine approximations, every real number $\alpha$ has a rational approximation $a/q$, where $a$ and $q$ are integers subject to
\begin{equation}\label{iii.2}
  1 \le q \le Q, \quad (a,q) = 1, \quad |\alpha - a/q| < (qQ)^{-1}.
\end{equation}

\begin{lemma}\label{l3.1}
  Let $k \ge 3$ and $0 < \rho < \sigma_k/(2+2\sigma_k)$. Suppose that $\alpha$ is real and that there exist integers $a$ and $q$ such that \eqref{iii.2} holds with $Q$ given by \eqref{iii.1}. 
  Let $|\xi_m| \le 1$, $|\eta_n| \le 1$, and define
  \[
    S(\alpha) = \sum_{m \sim M} \sum_{mn \in \mathcal I}
    \xi_m \eta_n e \left( \alpha (mn)^k \right).
  \]
  Then
  \[
    S(\alpha) \ll x^{\theta - \rho + \eps} + 
    \frac {w_k(q)^{1/2}x^{\theta + \eps}}{\left( 1 + \delta^2 x^k|\alpha - a/q| \right)^{1/2}},
  \]
  provided that
  \begin{equation}\label{iii.3}
    \delta^{-1} \max \left( x^{2\rho/\sigma_k}, \delta^{-k}x^{4\rho}, \big( \delta^{2k-2} x^{k-1+4k\rho} \big)^{1/(2k - 1)} \right) \ll M \ll x^{\theta-2\rho}.
  \end{equation}
\end{lemma}

\begin{proof}
  Set $H = \delta M$ and $N = xM^{-1}$ and define $\nu$ by $H^{\nu} = x^{2\rho}L^{-1}$. By~\eqref{iii.3}, we have $\nu < \sigma_k$. For $n_1, n_2 \le 2N$, let 
  \[
    \mathcal M(n_1, n_2) = \big\{ m \in (M, 2M] : mn_1, mn_2 \in \mathcal I \big\}.
  \]
  By Cauchy's inequality and an interchange of the order of summation,
  \begin{align}\label{iii.4}
    |S(\alpha)|^2  \ll x^{\theta}M + MT_1(\alpha),
  \end{align}
  where
  \[
    T_1(\alpha) = \sum_{n_1 < n_2}
    \left| \sum_{m \in \mathcal M(n_1, n_2)} e\left( \alpha \left( n_2^k - n_1^k \right) m^k \right) \right|.
  \]
  Let $\mathcal N$ denote the set of pairs $(n_1, n_2)$ with $n_1 < n_2$ and $\mathcal M(n_1, n_2) \neq \varnothing$ for which there exist integers $b$ and $r$ such that
  \begin{equation}\label{iii.5}
    1 \le r \le H^{k\nu}, \quad (b, r) = 1, \quad
    \left| r\left( n_2^k - n_1^k \right)\alpha - b \right| \le H^{k\nu}(\delta M^k)^{-1}.
  \end{equation}
  We remark that $\mathcal N$ contains $O(\delta N^2)$ pairs $(n_1,n_2)$. Since $\nu < \sigma_k$ and $M^k \le H^{k+1-2\nu}$, we can apply Lemma \ref{l2} with $\rho = \nu$, $x = M$ and $y = H$ to the inner summation in $T_1(\alpha)$. We get
  \begin{equation}\label{iii.6}
    T_1( \alpha) \ll x^{2\theta - 2\rho + \eps}M^{-1} + T_2(\alpha),
  \end{equation}
  where
  \[
    T_2( \alpha) = \sum_{(n_1, n_2) \in \mathcal N}
    \frac {w_k(r)H}{1 + \delta M^k \left| \left( n_2^k - n_1^k \right)\alpha - b/r \right|}.
  \]

  We now change the summation variables in $T_2(\alpha)$ to
  \[
    d = (n_1, n_2), \quad n = n_1/d, \quad h = (n_2 - n_1)/d.
  \]
  We obtain
  \begin{equation}\label{iii.7}
    T_2( \alpha) \ll \sum_{dh \le \delta N} \sum_n \frac{ w_k(r)H }
    {1 + \delta M^k \left| hd^k R(n, h) \alpha - b/r \right|},
  \end{equation}
  where $R(n, h) = \left( (n + h)^k - n^k \right) / h$ and the inner summation is over $n$ with $(n, h) = 1$ and $(nd, (n+h)d) \in \mathcal N$. For each pair $(d, h)$ appearing in the summation on the right side of \eqref{iii.7}, Dirichlet's theorem on Diophantine approximation yields integers $b_1$ and $r_1$ with
  \begin{equation}\label{iii.8}
    1 \le r_1 \le x^{-2k\rho}(\delta M^{k}), \quad (b_1, r_1) = 1, \quad
    \left| r_1hd^k \alpha - b_1 \right| \le x^{2k\rho}(\delta M^k)^{-1}.
  \end{equation}
  As $R(n, h) \le 3^kN^{k - 1}$, combining \eqref{iii.3}, \eqref{iii.5} and \eqref{iii.8}, we get
  \begin{align*}
    \left| b_1r R(n, h) - br_1 \right|
    & \le r_1 H^{k\nu}(\delta M^k)^{-1} + r R(n, h)x^{2k\rho}(\delta M^k)^{-1} \\
    & \le L^{-k} + 3^{k}\delta^{-1}x^{k-1+4k\rho}M^{1-2k}L^{-k} < 1.
  \end{align*}
  Hence,
  \begin{equation}\label{iii.9}
    \frac br = \frac {b_1 R(n, h)}{r_1}, \quad r = \frac {r_1}{(r_1, R(n, h))}.
  \end{equation}
  Combining \eqref{iii.7} and \eqref{iii.9}, we obtain
  \[
    T_2( \alpha) \ll \sum_{dh \le \delta N}
    \frac H{1 + \delta M^kN_d^{k - 1} \left| hd^k \alpha - b_1/r_1 \right|}
    \sum_{ \substack{ n \sim N_d\\ (n, h) = 1}}
    w_k \left( \frac {r_1}{(r_1, R(n, h))} \right),
  \]
  where $N_d = Nd^{-1}$. Using \eqref{ii.4}, we deduce that
  \begin{equation}\label{iii.10}
    T_2( \alpha) \ll \delta x^{\theta + \eps} + T_3(\alpha),
  \end{equation}
  where
  \[
    T_3(\alpha) = \sum_{dh \le \delta N}
    \frac{ r_1^{\eps}w_k(r_1)HN_d }{1 + \delta M^kN_d^{k - 1} \left| hd^k \alpha - b_1/r_1 \right|}.
  \]

  We now write $\mathcal H$ for the set of pairs $(d, h)$ with $dh \le \delta N$ for which there exist integers $b_1$ and $r_1$ subject to
  \begin{equation}\label{iii.11}
    1 \le r_1 \le x^{2k\rho}, \quad (b_1, r_1) = 1, \quad
    \left| r_1 hd^k \alpha - b_1 \right| \le x^{-k+1+2k\rho}H^{-1}.
  \end{equation}
  We have
  \begin{equation}\label{iii.12}
    T_3( \alpha) \ll  x^{2\theta - 2\rho + \eps}M^{-1} + T_4(\alpha),
  \end{equation}
  where
  \[
    T_4(\alpha) = \sum_{(d, h) \in \mathcal H}
    \frac{ r_1^{\eps}w_k(r_1)HN_d }{1 + \delta M^kN_d^{k - 1} \left| hd^k \alpha - b_1/r_1 \right|}.
  \]
  For each $d \le \delta N$, Dirichlet's theorem on Diophantine approximation yields integers $b_2$ and $r_2$ with
  \begin{equation}\label{iii.13}
    1 \le r_2 \le {\textstyle \frac 12}x^{k-1-2k\rho}H, \quad (b_2, r_2) = 1, \quad
    \left| r_2d^k \alpha - b_2 \right| \le 2x^{-k+1+2k\rho}H^{-1}.
  \end{equation}
  Combining \eqref{iii.11} and \eqref{iii.13}, we obtain
  \begin{align*}
    \left| b_2r_1h - b_1r_2 \right| &\le
      (r_2+2r_1h)x^{-k+1+2k\rho}H^{-1} \\
    & \le {\textstyle \frac 12} + 2x^{-k+2+4k\rho}M^{-2} < 1,
  \end{align*}
  whence
  \[
    \frac {b_1}{r_1} = \frac {h b_2}{r_2}, \quad
    r_1 = \frac {r_2}{(r_2, h)}.
  \]
  We write $Z_d = \delta M^kN_d^{k - 1} \left| d^k \alpha - b_2 / r_2 \right|$ and we use \eqref{ii.3} to get
  \[
    T_4(\alpha) \le \sum_{dh \le \delta N} 
    \frac { r_2^{\eps}HN_d }{1 + Z_dh} w_k \left( \frac {r_2}{(r_2, h)} \right) 
    \ll \sum_{d \le \delta N} \frac { w_k(r_2) x^{2\theta+\eps}M^{-1} }{d^2(1 + \delta Z_dN_d)}.
  \]
  Hence,
  \begin{equation}\label{iii.14}
    T_4(\alpha) \ll x^{2\theta - 2\rho + \eps}M^{-1} + T_5(\alpha),
  \end{equation}
  where
  \[
    T_5(\alpha) = \sum_{d \in \mathcal D}
    \frac{ w_k(r_2) x^{2\theta + \eps}M^{-1} }{d^2 \left( 1 + \delta^2(x/d)^k
    \left| d^k\alpha - b_2/r_2 \right| \right)}
  \]
  and $\mathcal D$ is the set of integers $d \le x^{2\rho}$ for which there exist integers $b_2$ and $r_2$ with
  \begin{equation}\label{iii.15}
    1 \le r_2 \le x^{2k\rho}, \quad (b_2, r_2) = 1, \quad
    \left| r_2d^k\alpha - b_2 \right| \le \delta^{-2}x^{-k+2k\rho}.
  \end{equation}
  Combining \eqref{iii.1}, \eqref{iii.2} and \eqref{iii.15}, we deduce that
  \begin{align*}
    \left| r_2 d^k a - b_2q \right| &\le r_2d^kQ^{-1} + q\delta^{-2}x^{-k+2k\rho} \\
    & \le x^{4k\rho}Q^{-1} + \delta^{-2}x^{-k+2k\rho}Q < 1,
  \end{align*}
  whence
  \[
    \frac {b_2}{r_2} = \frac {d^ka}q, \quad r_2 = \frac q{(q, d^k)}.
  \]
  Thus, recalling \eqref{ii.3}, we get
  \begin{equation}\label{iii.16}
    T_5(\alpha) \ll \frac{ x^{2\theta + \eps}M^{-1} }
    {1 + \delta^2x^{k}\left| \alpha - a/q \right|} \sum_{d \le x^{2\rho}}
    w_k \left( q / (q, d^k) \right) d^{-2} \ll
    \frac{ w_k(q)x^{2\theta + \eps}M^{-1} } {1 + \delta^2x^{k}\left| \alpha - a/q \right|}.
  \end{equation}

  The lemma follows from \eqref{iii.3}, \eqref{iii.4}, \eqref{iii.6},
  \eqref{iii.10}, \eqref{iii.12}, \eqref{iii.14} and \eqref{iii.16}.
\end{proof}

\begin{lemma}\label{l3.2}
  Let $k \ge 3$ and $0 < \rho < \sigma_k$. Suppose that $\alpha$ is real and that there exist integers $a$ and $q$ such that \eqref{iii.2} holds with $Q$ given by \eqref{iii.1}. Let $|\xi_{m_1,m_2}| \le 1$, and define
  \[
    S(\alpha) = \sum_{m_1 \sim M_1} \sum_{m_2 \sim M_2} \sum_{m_1m_2n \in \mathcal I} \xi_{m_1,m_2}
    e \left( \alpha (m_1m_2n)^k \right).
  \]
  Then
  \[
    S(\alpha) \ll x^{\theta - \rho + \eps}
    + \frac {w_k(q)x^{\theta + \eps}}{1 + \delta x^{k}|\alpha - a/q|},
  \]
  provided that
  \begin{equation}\label{iii.17}
    M_1^{2k-1} \ll \delta x^{k-(2k+1)\rho}, \quad
    M_1M_2 \ll \min(\delta x^{1 - \rho/\sigma_k}, \delta^{k+1}x^{1-2\rho}), \quad
    M_1M_2^2 \ll \delta^{1/k}x^{1 - 2\rho}.
  \end{equation}
\end{lemma}

\begin{proof}
  Set $N = x(M_1M_2)^{-1}$ and $H = \delta N$ and  define $\nu$ by $H^{\nu} = x^{\rho}L^{-1}$. Note that, by \eqref{iii.17}, we have $\nu < \sigma_k$. We denote by $\mathcal M$ the set of pairs $(m_1, m_2)$, with $m_1 \sim M_1$ and $m_2 \sim M_2$, for which there exist integers $b_1$ and $r_1$ with
  \begin{equation}\label{iii.18}
    1 \le r_1 \le H^{k\nu}, \quad (b_1, r_1) = 1, \quad
    \left| r_1 (m_1m_2)^k \alpha - b_1 \right| \le H^{k\nu}(\delta N^k)^{-1}.
  \end{equation}
  We apply Lemma \ref{l2} to the summation over $n$ and get
  \begin{equation}\label{iii.19}
    S(\alpha) \ll x^{\theta - \rho + \eps} + T_1(\alpha),
  \end{equation}
  where
  \[
    T_1(\alpha) = \sum_{(m_1, m_2) \in \mathcal M} \frac {w_k(r_1)H}
    {1 + \delta N^{k} \left| (m_1m_2)^k \alpha - b_1/r_1 \right|}.
  \]
  For each $m_1 \sim M_1$, we apply Dirichlet's theorem on Diophantine approximation to find integers $b$ and $r$ with
  \begin{equation}\label{iii.20}
    1 \le r \le x^{-k\rho}(\delta N^k), \quad (b, r) = 1, \quad
    \left| r m_1^k \alpha - b \right| \le x^{k\rho}(\delta N^k)^{-1}.
  \end{equation}
  By \eqref{iii.17}, \eqref{iii.18} and \eqref{iii.20},
  \begin{align*}
    \left| b_1r - bm_2^kr_1 \right| & \le rH^{k\nu}(\delta N^k)^{-1} + r_1m_2^kx^{k\rho}(\delta N^k)^{-1} \\
    & \le L^{-k} + 2^{k} \delta^{-1}x^{-k + 2k\rho}(M_1M_2^2)^kL^{-k} < 1,
  \end{align*}
  whence
  \[
    \frac {b_1}{r_1} = \frac {m_2^kb}{r}, \quad r_1 = \frac r{(r, m_2^k)}.
  \]
  Thus, by \eqref{ii.3},
  \begin{align}\label{iii.21}
    T_1(\alpha) & \ll \sum_{m_1 \sim M_1} \frac {H}{1 + \delta (M_2N)^k\left| m_1^k \alpha - b/r \right|}
    \sum_{m_2 \sim M_2} w_k \left( \frac r{(r, m_2^k)} \right)\\
    &\ll \sum_{m_1 \sim M_1} \frac {r^{\eps}w_k(r)HM_2}
    {1 + \delta (M_2N)^k \left| m_1^k \alpha - b/r \right|}. \notag
  \end{align}
  Let $\mathcal M_1$ be the set of integers $m \sim M_1$ for which there exist integers $b$ and $r$ with
  \begin{equation}\label{iii.22}
    1 \le r \le x^{k\rho}L^{-1}, \quad (b, r) = 1, \quad
    \left| r m^k \alpha - b \right| \le \delta^{-1} x^{-k+k\rho}M_1^kL^{-1}.
  \end{equation}
  From \eqref{iii.21},
  \begin{equation}\label{iii.23}
    T_1(\alpha) \ll x^{\theta - \rho + \eps} + T_2(\alpha),
  \end{equation}
  where
  \[
    T_2(\alpha) = \sum_{m \in \mathcal M_1} 
    \frac {r^{\eps}w_k(r)HM_2}{1 + \delta (M_2N)^k \left| m^k \alpha - b/r \right|}.
  \]
  We now consider two cases depending on the size of $q$ in \eqref{iii.2}.

  \medskip

  \noindent
  {\em Case 1:} $q \le \delta x^{k - k\rho}M_1^{-k}$. In this case, we estimate $T_2(\alpha)$ as in the proof of Lemma~\ref{l3.1}. Combining \eqref{iii.1}, \eqref{iii.2}, \eqref{iii.17} and \eqref{iii.22}, we obtain
  \begin{align*}
    \left| rm^ka - bq \right| &\le q\delta^{-1} x^{-k+k\rho}M_1^kL^{-1} + rm^kQ^{-1} \\
    &\le L^{-1} + 2^k x^{k\rho}M_1^kQ^{-1}L^{-1} < 1.
  \end{align*}
  Therefore,
  \[
    \frac br = \frac {m^ka}q, \quad r = \frac q{(q, m^k)},
  \]
  and by \eqref{ii.3},
  \begin{equation}\label{iii.24}
    T_2(\alpha) \ll \frac {q^{\eps}HM_2}{1 + \delta x^k |\alpha - a/q|}
    \sum_{m \sim M_1} w_k \left( \frac q{(q, m^k)} \right)
    \ll \frac {w_k(q)x^{\theta + \eps}}{1 + \delta x^k |\alpha - a/q|}.
  \end{equation}

  \medskip

  \noindent
  {\em Case 2:} $q > \delta x^{k - k\rho}M_1^{-k}$. We remark that in this case, the choice \eqref{iii.1} implies that $M_1 \ge x^{\rho}$. By a standard splitting argument,
  \begin{equation}\label{iii.25}
    T_2(\alpha) \ll \sum_{d \mid q} \sum_{m \in \mathcal M_d(R, Z)}
    \frac {w_k(r)HM_2x^{\eps}}{1 + \delta(M_2N)^k(RZ)^{-1}},
  \end{equation}
  where
  \begin{equation}\label{iii.26}
    1 \le R \le x^{k\rho}L^{-1}, \quad \delta x^{k - k\rho}M_1^{-k}L \le Z \le \delta (x/M_1)^k,
  \end{equation}
  and $\mathcal M_d(R, Z)$ is the subset of $\mathcal M_1$ containing integers $m$ subject to
  \[
    (m, q) = d, \quad r \sim R, \quad \left| rm^k\alpha - b \right| < Z^{-1}.
  \]

  We now estimate the inner sum on the right side of \eqref{iii.25}. We have
  \begin{equation}\label{iii.27}
    \sum_{m \in \mathcal M_d(R, Z)} w_k(r) \ll \sum_{r \sim R} w_k(r) \mathcal S_0 (r),
  \end{equation}
  where $\mathcal S_0 (r)$ is the number of integers $m \sim M_1$ with $(m, q) = d$ for which there exists an integer $b$ such that
  \begin{equation}\label{iii.28}
    (b, r) = 1 \quad \text{and} \quad \left| r m^k \alpha - b \right| < Z^{-1}.
  \end{equation}
  Since for each $m \sim M_1$ there is at most one pair $(b, r)$ satisfying \eqref{iii.28} and $r \sim R$, we have
  \begin{equation}\label{iii.29}
    \sum_{r \sim R} \mathcal S_0(r) \le \sum_{ \substack{ m \sim M_1\\ (m, q) = d}} 1
    \ll M_1d^{-1} + 1.
  \end{equation}
  Hence,
  \begin{align}\label{iii.30}
    \sum_{ \substack{ r \sim R\\ (q, rd^k) = q}} w_k(r)\mathcal S_0(r)
    \ll R^{-1/k} \left( M_1d^{-1} + 1 \right) \ll M_1q^{-1/k} + 1,
  \end{align}
  on noting that the sum on the left side is empty unless $Rd^k \gg q$. 
  
  When $(q, rd^k) < q$, we make use of Lemma \ref{l3}. By \eqref{iii.2}, \eqref{iii.26} and \eqref{iii.28},
  \begin{equation}\label{iii.31}
    \mathcal S_0(r) \le \mathcal S(r),
  \end{equation}
  where we $\mathcal S(r)$ is the number of integers $m$ subject to
  \[
    m \sim M_1d^{-1}, \quad (m, q_1) = 1, \quad \left\| ard^{k - 1}m^k / q_1 \right\| < \Delta,
  \]
  with $q_1 = qd^{-1}$ and $\Delta = Z^{-1} + 2^{k+1}RM_1^k(qQ)^{-1}$. Since \eqref{iii.17} implies $M_1 \le \delta x^{k - k\rho}M_1^{-k} < q$, we obtain
  \begin{equation}\label{iii.32}
    \mathcal S(r) \ll \Delta q^{\eps}d^{-1} (M_1 + q) \ll \Delta q^{1 + \eps}.
  \end{equation}
  Combining \eqref{iii.31} and \eqref{iii.32}, we get
  \begin{equation}\label{iii.33}
    \mathcal S_0(r) \ll \Delta q^{1 + \eps}.
  \end{equation}
  We now apply H\"older's inequality, \eqref{ii.2}, \eqref{iii.29}, and \eqref{iii.33} and obtain
  \begin{align}\label{iii.34}
    \sum_{ \substack{ r \sim R\\ (q, rd^3) < q}} w_3(r)\mathcal S_0(r)
      & \ll \left( \Delta q^{1 + \eps} \right)^{1/4}
      \left( \sum_{r \sim R} w_3(r)^4 \right)^{1/4}
      \left( \sum_{r \sim R} \mathcal S_0(r) \right)^{3/4} \\
    & \ll \Delta^{1/4}q^{1/4 + \eps} R^{-1/4} M_1^{3/4}.
      \notag
  \end{align}
  Similarly, when $k \ge 4$, we have
  \begin{align}\label{iii.35}
    \sum_{ \substack{ r \sim R\\ (q, rd^k) < q}} w_k(r)\mathcal S_0(r)
      & \ll \left( \Delta q^{1 + \eps} \right)^{1/k}
      \left( \sum_{r \sim R} w_k(r)^k \right)^{1/k}
      \left( \sum_{r \sim R} \mathcal S_0(r) \right)^{1 - 1/k} \\
    & \ll \Delta^{1/k}q^{1/k + \eps} R^{(1 - k)/k^2} M_1^{(k - 1)/k}.
      \notag
  \end{align}
  Combining \eqref{iii.27}, \eqref{iii.30}, \eqref{iii.34} and \eqref{iii.35}, we deduce
  \begin{equation}\label{iii.36}
    \sum_{m \in \mathcal M_d(R, Z)} w_3(r) \ll
    \Delta^{1/4}q^{1/4 + \eps} R^{-1/4} M_1^{3/4} + M_1q^{-1/3} + 1
  \end{equation}
  and
  \begin{equation}\label{iii.37}
    \sum_{m \in \mathcal M_d(R, Z)} w_k(r) \ll
    \Delta^{1/k}q^{1/k + \eps} R^{(1 - k)/k^2} M_1^{(k - 1)/k}
    + M_1q^{-1/k} + 1
  \end{equation}
  for $k \ge 4$.

  Substituting \eqref{iii.36} into \eqref{iii.25}, we get
  \begin{align*}
    T_2(\alpha) &\ll \frac {x^{\theta+\eps}M_1^{-1/4}}
      {1 + \delta (M_2N)^3(RZ)^{-1}}
      \left( \frac Q{RZ} + \frac {M_1^3}Q \right)^{1/4} +
      x^{\theta + \eps}q^{-1/3} + x^{\theta + \eps}M_1^{-1} \\
    &\ll (\delta^3 xM_1^2Q)^{1/4 + \eps}
      + x^{\theta + \eps} \left( M_1^2Q^{-1} \right)^{1/4}
      + x^{\rho + \eps}M_1 + x^{\theta - \rho + \eps}. 
  \end{align*}
  The hypotheses of the lemma ensure that 
  \[
    M_1 \le \min\left( \delta^{1/2}x^{3/2-2\rho}Q^{-1/2}, Q^{1/2}x^{-2\rho}, x^{\theta-2\rho} \right),
  \]
  and so when $k = 3$,
  \begin{equation}\label{iii.38}
    T_2(\alpha) \ll x^{\theta - \rho + \eps}.
  \end{equation}
  When $k \ge 4$, by \eqref{iii.25} and \eqref{iii.37},
  \begin{align*}
    T_2(\alpha) &\ll \frac {x^{\theta+\eps}M_1^{-1/k}R^{1/k^2}} {1 + \delta (M_2N)^k(RZ)^{-1}}
      \left( \frac Q{RZ} + \frac {M_1^k}Q \right)^{1/k} +
      x^{\theta + \eps}q^{-1/k} + x^{\theta + \eps}M_1^{-1} \\
    &\ll \left( x^{\rho}Q(\delta M_1)^{k - 1} \right)^{1/k + \eps}
      + x^{\theta + \eps} \left( x^{\rho}M_1^{k - 1}Q^{-1} \right)^{1/k}
      + x^{\rho + \eps}M_1 + x^{\theta - \rho + \eps},
  \end{align*}
  and using \eqref{iii.1} and \eqref{iii.17}, we find that \eqref{iii.38} holds in this case as well.

  The desired estimate follows from \eqref{iii.19}, \eqref{iii.23}, \eqref{iii.24} and \eqref{iii.38}.
\end{proof}

\section{Proof of Theorems \ref{th2} and \ref{th3}}

In this section we deduce the main theorems from Lemmas \ref{l3.1} and \ref{l3.2} and Heath-Brown's identity for $\Lambda(n)$. We apply Heath-Brown's identity in the following form \cite[Lemma~1]{HB82}: if $n \le X$ and $J$ is a positive integer, then
\begin{equation}\label{iv.1}
  \Lambda(n) = \sum_{j = 1}^{J} \binom {J}j (-1)^j 
  \sum_{\substack{ n = n_1 \cdots n_{2j}\\ n_1, \dots, n_j \le X^{1/J}}} \mu(n_1) \cdots \mu(n_j) (\log n_{2j}).
\end{equation}

Let $\alpha \in \mathfrak m(P)$. By Dirichlet's theorem on Diophantine approximation, there exist integers $a$ and $q$ such that \eqref{iii.2} holds with $Q$ given by \eqref{iii.1}. Let $\beta$ be defined by
\[
  x^{\beta} = \min\left( \delta^2 x^{1 - 2\rho(\sigma_k^{-1}+1)}, \delta^{k+2}x^{1-6\rho},
  \big( \delta^{2k}x^{k-(8k-2)\rho} \big)^{1/(2k-1)} \right),
\]
and suppose that $\rho$ and $\delta$ are chosen so that
\begin{equation}\label{iv.2}
  \delta^{-1}x^{\beta + 2\rho} \ge 2x^{1/3}.
\end{equation}
We apply \eqref{iv.1} with $X = x+x^{\theta}$ and $J \ge 3$ chosen so that $x^{1/J} \le x^{\beta}$. After a standard splitting argument, we have
\begin{equation}\label{iv.3}
  \sum_{n \in \mathcal I} \Lambda(n)e\left( \alpha n^k \right) 
  \ll \sum_{\mathbf N} \left| \sum_{n \in \mathcal I} c(n; \mathbf N) e\left( \alpha n^k \right) \right|,
\end{equation}
where $\mathbf N$ runs over $O( L^{2J-1} )$ vectors $\mathbf N = (N_1, \dots, N_{2j})$, $j \le J$, subject to
\[
  N_1, \dots, N_j \ll x^{1/J}, \qquad x \ll N_1 \cdots N_{2j} \ll x,
\]
and 
\[
  c(n; \mathbf N) = \sum_{\substack{ n = n_1 \cdots n_{2j}\\ n_i \sim N_i}} \mu(n_1) \cdots \mu(n_j)(\log n_{2j}).
\]
In fact, since the coefficient $\log n_{2j}$ can be removed by partial summation, we may assume that
\[
  c(n; \mathbf N) = L \sum_{\substack{ n = n_1 \cdots n_{2j}\\ N_i < n_i \le N_i'}} \mu(n_1) \cdots \mu(n_j),
\]
where $N_i < N_i' \le 2N_i$ (in reality, $N_i' = 2N_i$ except for $i = 2j$). We also assume (as we may) that the summation variables $n_{j+1}, \dots, n_{2j}$ are labeled so that $N_{j + 1} \le \cdots \le N_{2j}$. Next, we show that each of the sums occurring on the right side of \eqref{iv.3} satisfies the bound
\begin{equation}\label{iv.4}
  \sum_{n \in \mathcal I} c(n; \mathbf N) e\left( \alpha n^k \right) \ll
  x^{\theta-\rho+\eps} + \frac {w_k(q)^{1/2}x^{\theta+\eps}}{\left( 1 + \delta^2 x^k|\alpha - a/q| \right)^{1/2}}.
\end{equation}
The analysis involves several cases depending on the sizes of $N_1, \dots, N_{2j}$.

\medskip

\paragraph{\em Case 1:} $N_1 \cdots N_{j} \gg \delta^{-1}x^{2\rho}$. Since none of the $N_i$'s exceeds $x^{\beta}$, there must be a set of indices $S \subset \{1, \dots, j\}$ such that 
\begin{equation}\label{iv.5}
  \delta^{-1}x^{2\rho} \le \prod_{i \in S} N_i \le \delta^{-1}x^{\beta + 2\rho}. 
\end{equation}
Hence, we can rewrite $c(n; \mathbf N)$ in the form
\begin{equation}\label{iv.6}
  c(n; \mathbf N) = \sum_{\substack{ mr = n\\ m \asymp M}} \xi_m\eta_r,
\end{equation}
where $|\xi_m| \le \tau(m)^c$, $|\eta_r| \le \tau(r)^c$, and $M = \prod_{i \notin S} N_i$. By \eqref{iv.5}, $M$ satisfies \eqref{iii.3}, so \eqref{iv.4} follows from Lemma \ref{l3.1}.

\smallskip

\paragraph{\em Case 2:} $N_1 \cdots N_{j} < \delta^{-1}x^{2\rho}$, $j \le 2$. When $j = 1$, \eqref{iv.4} follows from Lemma \ref{l3.2} with $M_1 = N_1$, $M_2 = 1$ and $N = N_2$. When $j = 2$, we have
\begin{gather*}
  N_3 \le (x/N_1N_2)^{1/2} \le x^{1/2}, \quad N_1N_2N_3 \le (xN_1N_2)^{1/2} \le \delta^{-1}x^{1/2 + \rho}, \\
  (N_1N_2)^2N_3 \le x^{1/2}(N_1N_2)^{3/2} \le \delta^{-2}x^{1/2 + 3\rho}.
\end{gather*}
Hence, we can deduce \eqref{iv.4} from Lemma \ref{l3.2} with $M_1 = N_3$, $M_2 = N_1N_2$ and $N = N_4$, provided that
\begin{gather}
  x^{k-1/2} \le \delta x^{k-(2k+1)\rho}, \quad \delta^{-2}x^{1/2+3\rho} \le \delta^{1/k}x^{1-2\rho}, \label{iv.7}\\
  \delta^{-1} x^{1/2+\rho} \le \delta\min\left( x^{1-\rho/\sigma_k}, \delta^k x^{1-2\rho} \right). \label{iv.8}
\end{gather}

\smallskip

\paragraph{\em Case 3:} $N_1 \cdots N_{j} < \delta^{-1}x^{2\rho}$, $j \ge 3$. In this case, we have 
\[
  N_{j+1}, \dots, N_{2j-2} \le 2x^{1/3} \le \delta^{-1}x^{\beta + 2\rho}.
\]

\medskip

\paragraph{\em Case 3.1:} $N_1 \cdots N_{2j-2} \ge \delta^{-1}x^{2\rho}$. Let $r$ be the least index with $N_1 \cdots N_r \ge \delta^{-1}x^{2\rho}$. We can use the product $N_1 \cdots N_r$ in a similar fashion to the product $N_1 \cdots N_j$ in Case 1 to represent $c(n; \mathbf N)$ in the form \eqref{iv.6}. Thus, we can appeal to Lemma \ref{l3.1} to show that \eqref{iv.4} holds again. 

\medskip

\paragraph{\em Case 3.2:} $N_1 \cdots N_{2j-2} < \delta^{-1}x^{2\rho}$. Then we are in a similar situation to Case 2 with $j=2$, with the product $N_1 \cdots N_{2j-2}$ playing the role of $N_1N_2$ in Case 2. Thus, we can again use Lemma \ref{l3.2} to obtain \eqref{iv.4}. 

\medskip

By the above analysis, 
\begin{equation}\label{iv.9}
  \sum_{n \in \mathcal I} \Lambda(n)e\left( \alpha n^k \right) \ll
  x^{\theta-\rho+\eps} + \frac {w_k(q)^{1/2}x^{\theta+\eps}}{\left( 1 + \delta^2 x^k|\alpha - a/q| \right)^{1/2}},
\end{equation}
provided that conditions \eqref{iv.2}, \eqref{iv.7} and \eqref{iv.8} hold. Altogether, those conditions are equivalent to the inequality
\begin{align*}
  x^{\rho} \ll \min\left( (\delta^3x^2)^{\sigma_k/6}, (\delta^2 x)^{1/(4k+2)}, (\delta^2 x)^{\sigma_k/(1+\sigma_k)}, \delta^{(k+2)/6}x^{1/6}, \right. \\
  \left. \delta^{(k+1)/4}x^{1/6}, \delta^{(2k+1)/5k}x^{1/10}, \delta^{1/(4k)}x^{(k+1)/(12k)} \right). 
\end{align*}
We have
\begin{gather*}
  \delta^{(k+2)/6}x^{1/6} \le \delta^{(k+1)/4}x^{1/6}, \quad 
  (\delta^2 x)^{1/(4k+2)} \le \delta^{1/(4k)}x^{(k+1)/(12k)}, \\
  (\delta^3 x^2)^{\sigma_k/6} \le (\delta^2 x)^{\sigma_k/(1+\sigma_k)} \qquad \text{when } \delta \ge x^{-1/3},
\end{gather*}
so the third, fifth and seventh terms in the above minimum are superfluous. Recalling the definition of $\delta$, we conclude that \eqref{iv.9} holds whenever
\begin{equation}\label{iv.10}
  \rho \le  \min\left( \frac {\sigma_k(3\theta-1)}6, \frac {2\theta - 1}{4k+2}, 
  \frac{(k+2)\theta - k - 1}6, \frac {(4k+2)\theta - 3k-2}{10k} \right).
\end{equation}
The latter minimum is exactly the function $\rho_k(\theta)$ defined in the statements of Theorems \ref{th2} and \ref{th3}. Indeed, when $k = 3$, the first term in the minimum is always larger than the second, so it can be discarded and we are left with $\rho_3(\theta)$. On the other hand, when $k \ge 4$, the second and fourth terms in the minimum are superfluous. Therefore, \eqref{iv.10} is a direct consequence of the hypotheses of the theorems and the proof of \eqref{iv.9} is complete.

If either $q \ge x^{2k\rho}$ or $|q\alpha - a| \ge \delta^{-2}x^{k-2k\rho}$, we can use \eqref{ii.1} to show that the second term on the right side of \eqref{iv.9} is smaller than the first. Thus, 
\begin{equation}\label{iv.11}
  \sup_{\alpha \in \mathfrak m(x^{2k\rho})} \left| f_k \left(\alpha; x, x^{\theta} \right) \right| \ll x^{\theta-\rho+\eps}.
\end{equation}
This establishes the theorems when $P \ge x^{2k\rho}$. When $P < x^{2k\rho}$, Theorem \ref{th1} gives 
\[
  \sup_{\alpha \in \mathfrak m(P) \cap \mathfrak M(x^{2k\rho})} \left| f_k \left(\alpha; x, x^{\theta} \right) \right| \ll x^{\theta-\rho+\eps} + x^{\theta+\eps}P^{-1/2},
\]
which in combination with \eqref{iv.11} establishes the theorems in the case $P < x^{2k\rho}$. \qed

\end{document}